\documentclass[12pt]{amsart}
\usepackage[centertags]{amsmath}
\usepackage{amsfonts,amssymb}
\usepackage{graphicx}
\usepackage{enumerate}
\usepackage[latin1]{inputenc}
\usepackage{placeins}
\usepackage{float}

  \newtheorem{theorem}{Theorem}
  
  \newtheorem{proposition}{Proposition}
  \newtheorem{lemma}{Lemma}%
  \theoremstyle{remark}

  \newcommand{\ZZ}{\mathbb Z}
  \newcommand{\NN}{\mathbb N}
  \newcommand{\Sc}{\mathcal S}
   
\begin{document}

\title{Some thoughts on pseudoprimes}
\date{\today}

\author{Carl Pomerance}
\address{Mathematics Department, Dartmouth College, Hanover, NH 03784}
\email{carl.pomerance@dartmouth.edu}

\author{Samuel S. Wagstaff, Jr.}
\address{Center for Education and Research in Information Assurance and Security
  and Department of Computer Sciences, Purdue University \\
  West Lafayette, IN 47907-1398 USA}
\thanks{S.S.W.'s work was supported by the CERIAS Center at Purdue University}

\email{ssw@cerias.purdue.edu}

\begin{abstract}
We consider several problems about pseudoprimes.  First, we look at
the issue of their distribution in residue classes.  There is a literature
on this topic in the case that the residue class is coprime to the modulus.
Here we provide some robust statistics in both these cases and the general
case.  In particular we tabulate all even pseudoprimes to $10^{16}$.
Second, we prove a recent conjecture of Ordowski: the set of integers $n$
which are a pseudoprime to some base which is a proper divisor of $n$
has an asymptotic density.
\end{abstract} 
\maketitle
\vskip-30pt
\newenvironment{dedication}
        {\vspace{6ex}\begin{quotation}\begin{center}\begin{em}}
        {\par\end{em}\end{center}\end{quotation}}
\begin{dedication}
{In memory of Aleksandar Ivi\'c (1949--2020)}
\end{dedication}
\vskip20pt

\section{Introduction}

Fermat's ``little" theorem is part of the basic landscape in elementary number theory.
It asserts that if $p$ is a prime, then $a^p\equiv a\pmod p$ for every prime $p$.
One interest in this result is that for a given pair $a,p$, it is not hard computationally
to check if the congruence holds.  So, if the congruence fails, we have proved
that the modulus $p$ is not prime.

A pseudoprime is a composite number $n$ with $2^n\equiv2\pmod n$, and
more generally, a pseudoprime base $a$ is a composite number $n$ with
$a^n\equiv a\pmod n$.  Pseudoprimes exist, in fact, there are composite
numbers $n$ which are pseudoprimes to every base $a$, the first 3 examples
being 561, 1105, and 1729.   These are the
Carmichael numbers.  Named after Carmichael \cite{C} who published the first
few examples in 1910, they were actually anticipated by quite a few years by
\v{S}imerka \cite{S}.

We now know that there are infinitely many Carmichael numbers (see \cite{AGP}),
the number of them up to $x$ exceeding $x^{0.33}$ for all sufficiently large $x$
(see \cite{H}).  This count holds {\it a fortiori} for pseudoprimes to any
fixed base $a$ since the Carmichael numbers comprise a subset of the
base-$a$ pseudoprimes.

One can also ask for upper bounds on the distribution of pseudoprimes and
Carmichael numbers.  Let 
\[
L(x)=\exp(\log x\log\log\log x/\log\log x)=x^{\frac{\log\log\log x}{\log\log x}}.
\]
We know (see \cite{Ppsp})
that the number of Carmichael numbers up to $x$ is at most
$x/L(x)$ for all sufficiently large $x$, and it is conjectured that this is
almost best possible in that the count is of the form $x/L(x)^{1+o(1)}$
as $x\to\infty$.  The heuristic for this assertion is largely based on thoughts
of Erd\H os \cite{E}.

It is conjectured that the same is true for pseudoprimes to any fixed base $a$,
however the upper bound is not as tight.  We know (see \cite{Ppsp}) that
for all large $x$, the number of {\it odd} pseudoprimes up to $x$ is
$\le x/L(x)^{1/2}$ and it seems likely that the argument goes through for
those base-$a$ pseudoprimes coprime to $a$, for any $a>1$.   It is likely
too that similar methods could be used to get comparable upper bounds
where the coprime condition is relaxed.

For positive coprime integers $a,n$ let $l_a(n)$ denote the order of $a\pmod n$
in $(\ZZ/n\ZZ)^*$.  Further, let $\lambda(n)$ denote the maximal value of
$l_a(n)$ over all $a\pmod n$; it is the universal exponent for the group
$(\ZZ/n\ZZ)^*$.  If $a,n$ are positive integers, not necessarily coprime, let
$n_a$ denote the largest divisor of $n$ coprime to $a$.  Note that $n$ is
a base-$a$ pseudoprime if and only if $l_a(n_a)\mid n-1$ and $n/n_a\mid a$,
as is easily verified.

It is natural to consider the distribution of pseudoprimes in residue classes.
Consider the integers $n$ with $n\equiv r\pmod m$, and suppose that $n$ is a base-$a$
pseudoprime.  Let us write down some necessary conditions for this to
occur.  Let 
\[
g=\gcd(r,m),\quad h=\gcd(l_a(g_a),m).
\]
Then if $n$ is a base-$a$ pseudoprime in the residue class $r\pmod m$, 
we must have 
\begin{equation}
\label{eq:cond}
h\mid r-1~\hbox{ and }~ g/g_a\mid a. 
\end{equation}
Further, if $g$ is even, the residue class $r\pmod m$ must contain
some integer $k$ with the Jacobi symbol $(a/k_{2a})=1$.
 We conjecture that these conditions are
sufficient for there to be infinitely many base-$a$ pseudoprimes
$n\equiv r\pmod m$.
In fact, a heuristic argument based
on that of Erd\H os \cite{E} suggests that if these conditions hold for $a,r,m$, then
the number $P_{a,r,m}(x)$ of base-$a$ pseudoprimes $n\equiv r\pmod m$ with $n\le x$
is $x^{1-o(1)}$ as $x\to\infty$.   
We prove the necessity of these conditions in Section \ref{secres}.

Here are a few examples.  If $\gcd(r,m)=1$ then the conditions hold for any $a$,
and the conjecture asserts that there are infinitely many base-$a$ pseudoprimes
that are $\equiv r\pmod m$.  In fact, this is known, see below.  If $a=2,\,r=0,\,m=2$,
the conditions hold and we are looking at even (base-2) pseudoprimes.  The
first example was found by Lehmer, and their infinitude was proved by
Beeger, see below.   The criteria instantly
tell us there are no (base-2) pseudoprimes divisible by 4, since if $a=2,\,r=0,\,m=4$ we have $g=g/g_a=4$ and $4\,\nmid\,2$.
An interesting case is $a=2,\,r=15,\,m=20$.  Then $g=5$ and
$h=4$.  Since $4\,\nmid\,15-1$, the condition \eqref{eq:cond} fails, and indeed
there are no pseudoprimes in the class 15 (mod 20).  
Another interesting case
is $a=2,\,r=6,\,m=16$.  We have $g=2,\,g_a=1,\,h=1$, so that \eqref{eq:cond}
holds.  But any integer $k\equiv6\pmod{16}$ has $k_{2a}\equiv3\pmod8$,
and $(2/k_{2a})=-1$.  So, the final condition fails, and indeed there are
no base-2 pseudoprimes in the class $6\pmod{16}$.

Let $C_{r,m}(x)$ denote the number of Carmichael
numbers $n\le x$ with $n\equiv r\pmod m$.  Clearly for any $a,r,m$ we have
$C_{r,m}(x)\le P_{a,r,m}(x)$.

Here are some things we know towards the conjecture.
\begin{itemize}
\item
For all large $x$ we have $C_{0,1}(x)>x^{.33}$.  This is the main result
of Harman \cite{H}, improving the earlier result with exponent $2/7$ in
\cite{AGP}.  
\item
If $\gcd(r,m)=1$ and $r$ is a square mod $m$, then for $x$ sufficiently
large, $C_{r,m}(x)>x^{1/5}$.  This result is due to Matom\"aki \cite{M}.
\item
If $\gcd(r,m)=1$, then  $C_{r,m}(x)>x^{1/(6\log\log\log x)}$ for $x$ sufficiently large.
This recent result of the first-named author \cite{P} is based on the argument
for a somewhat
weaker bound due to Wright \cite{W1}.
\item
If $\gcd(r,m)=1$, then $P_{2,r,m}(x)$ is unbounded.  This result of
Rotkiewicz \cite{R} is, of course,
weaker than the previous item, but it preceded it by over half a century
and is much simpler.
\end{itemize}

There are elementary ideas for showing $P_{2,r,m}(x)$ is unbounded even
when $\gcd(r,m)>1$.  For example, there are infinitely many even pseudoprimes,
the case $r=0,m=2$.  Here's a proof.  Suppose $n$ is an even pseudoprime
and let $p$ be a prime with $l_2(p)=n$.  From Bang \cite{B} such a prime
$p$ exists.  Then $pn$ is another even pseudoprime.
It remains to note that $n=161{,}038$ is an even pseudoprime.  
This proof is
essentially due to Beeger \cite{Be}.  The example $161{,}038$ was found
by Lehmer in 1950.

A similar argument can be found for other choices of $r,m$, but we know
no general proof that $P_{a,r,m}(x)$ is unbounded when \eqref{eq:cond}
holds.

At the end of this paper we present substantial counts of pseudoprimes in
residue classes.

The usual thought with pseudoprimes is to fix the base $a$ and look at
pseudoprimes $n$ to the base $a$.  Instead, one can take the opposite
perspective and fix $n$, looking then at the bases $a$ for which $n$ is
a pseudoprime.  
 Let 
 \[
 F(n)=\#\{a\kern-5pt\pmod{n}:a^{n-1}\,\equiv\,1\kern-5pt\pmod{n}\}.
 \]
   From 
 Baillie--Wagstaff \cite{BW} and Monier \cite{Mo},
 we have
 \[
 F(n)=\prod_{p\mid n}\gcd(p-1,n-1),
 \]
 where $p$ runs over primes.
 Now let 
 \[
 F^*(n)=\#\{a\kern-5pt\pmod{n}:a^n\,\equiv\, a\kern-5pt\pmod{n}\}.
 \]
   Note that $F^*(n)=n$
 if and only if $n=1$, $n$ is a prime, or $n$ is a Carmichael number.
  The Baillie--Wagstaff and Monier
 formula can be enhanced as follows:
  \[
 F^*(n)=\prod_{p\mid n}(1+\gcd(p-1,n-1)).
 \]
 Note that $F^*(n)-F(n)$ is the number of residues $a\pmod n$ with
 $a^n\equiv a\pmod n$ and $\gcd(a,n)>1$.  Among these it is
 interesting to consider those $a$ that divide $n$.
 Let 
 \[
 D(n)=\#\{a\mid n:1<a<n,~a^n\,\equiv\, a\kern-5pt\pmod{n}\}
 \]
  and let
 \[
 \Sc=\{n\in\NN:\#D(n)>0\}.
 \]
 T. Ordowski \cite{O} has conjectured that $\Sc$ has an asymptotic density;
 counts up to $10^8$ by A. Eldar suggest that this density may be about $\frac58$.
 In Section \ref{seccomp} we present a proof that the density of $\Sc$ exists and
 is less than 1.
 
 \section{Proof of Ordowski's conjecture}
 
 For each integer $b\ge2$ let
 \[
 \Sc_b=\{ab:a\ge2,~a^{ab}\,\equiv\, a\kern-5pt\pmod{ab}\},
 \]
Then
\[
\Sc=\bigcup_{b\ge2}\Sc_b.
\]
Indeed, if $b\ge 2$ and $n=ab\in\Sc_b$, then $a\in D(n)$, so
$n\in\Sc$.  Conversely, if $n\in\Sc$ and $a\mid n$ with $1<a<n$ and
$a^n\equiv a\pmod n$, then $n\in\Sc_{n/a}$.  

We also remark that if $n=ab\in\Sc_b$, then $\gcd(b,a)=1$.  Indeed, if
$p$ is a common prime factor with $p^\alpha\,\|\,a$, then we have
$p^{\alpha+1}\mid n$ and $p^{\alpha+1}\mid a^n$, contradicting
$a^n\equiv a\pmod n$.

For a set $\mathcal A$ of positive integers, let 
$\delta(\mathcal A)$ be the asymptotic density of $\mathcal A$ should it exist.

\begin{proposition}
\label{prop:ssum}
For each integer $b\ge2$, $\delta(\Sc_b)$ exists and
\begin{equation}
\label{eq:Sbsum}
c_1:=\sum_{b\ge2}\delta(\Sc_b)<\infty.
\end{equation}
\end{proposition}
 \begin{proof}
 To see that $\delta(\Sc_b)$ exists we will show that
$\Sc_b\cup\{b\}$ is a finite union of residue classes.
 
   To get a feel for things, we work out the first few $b$'s.  The case $b=2$
 is particularly simple.  For $n$ to be in $\Sc_2$ it is necessary that
 $n/2$ be odd, since we need $\gcd(b,n/b)=1$.  And this condition is
 sufficient when $n>2$: it is easy to check that $(n/2)^n\equiv n/2\pmod n$.  Indeed
 the congruence is trivial modulo $n/2$ and it is trivial modulo 2.
 Thus $\Sc_2$ is the set of numbers that are $2\pmod 4$ (other than 2), with
 density $\frac14$.
 
 Now take $b=3$.  For $ab\in\Sc_3$ we consider the two cases $a\equiv1\pmod 3$,
 $a\equiv2\pmod 3$.  Every number of the form $3a$ with $a\equiv1\pmod3$ and
 $a>1$ is in $\Sc_3$, which gives density $\frac19$.  For $a\equiv2\pmod 3$ we need
$2^{3a}\equiv 2\pmod 3$ and this holds if and only if $a$ is odd.  That is,
$a\equiv5\pmod6$, and this condition is sufficient.  This part of $\Sc_3$ has
density $\frac1{18}$, so $\delta(\Sc_3)=\frac16$. 

We now work out the general structure of $\Sc_b$.  We have a number
$ab$, where $\gcd(a,b)=1$ and $a>1$.  We trivially have $a^{ab}\equiv a\pmod a$, so
the important condition is $a^{ab}\equiv a\pmod b$.  Since $\gcd(a,b)=1$,
this is equivalent to $a^{ab-1}\equiv 1\pmod b$, which holds if and only
if $d\mid ab-1$, where $d$ is the multiplicative order of $a\pmod b$.
This cannot hold unless $\gcd(d,b)=1$, and in this case, $a$ is in a residue class (mod~$d$).
So, if $a\equiv a_0\pmod b$ and $a_0\pmod b$ has multiplicative order $d$
with $\gcd(d,b)=1$, then such $a$'s lie in a residue class of modulus $bd$.
Thus, for each residue in $a_0\in(\ZZ/b\ZZ)^*$ with multiplicative order $d$
coprime to $b$ we have a residue class of modulus $b^2d$ consisting of 
all $ab\in\Sc_b$ with $a\equiv a_0\pmod b$ and $a\equiv b^{-1}\pmod d$.

Let $\lambda(b)$ denote the universal exponent for the group $(\ZZ/b\ZZ)^*$.
Thus, the divisors of $\lambda(b)$ run over all of the possible multiplicative
orders for elements in the group.
For $d\mid\lambda(b)$, let $N(d,b)$ denote the
number of elements $a_0\pmod b$ with multiplicative order $d$.   Thus,
\begin{equation}
\label{eq:sbdens}
\delta(\Sc_b)=\sum_{\substack{d\mid\lambda(b)\\\gcd(d,b)=1}}\frac{N(d,b)}{b^2d}.
\end{equation}

It seems difficult to work out a formula for $N(d,b)$ but we do have the relation
\begin{equation}
\label{eq:Ndb}
\sum_{d\mid\lambda(b)}N(d,b)=\varphi(b),
\end{equation}
which just reflects the partitioning of $(\ZZ/b\ZZ)^*$ by the orders of its elements.
We consider various cases.   First suppose that $\lambda(b)$
is smooth, more specifically, assume that $P(\lambda(b))<B(b):=\exp((\log b)^{1/2})$,
where $P(n)$ denotes the largest prime factor of $n$.  Note that the primes
dividing $\lambda(b)$ are the same primes that divide $\varphi(b)$, so
that $P(\varphi(b))<B(b)$.  Using the main result from
\cite{BFPS}, the number of such integers $b\le x$ is $\le x/B(x)$
for all sufficiently large $x$.  Since \eqref{eq:Ndb} implies that the sum of $N(d,b)/d$ for $d\mid\lambda(b)$
is $\le\varphi(b)<b$,  \eqref{eq:sbdens} implies that $\delta(\Sc_b)<1/b$.  But the
sum of $1/b$ over such a sparse set of $b$'s is easily seen to converge via
a partial summation argument.

So, we may assume that $p_b:=P(\lambda(b))\ge B(b)$.  There are two types of
numbers $d\mid\lambda(b)$ to consider: $p_b\mid d$ and $p_b\,\nmid\, d$.
In the first case \eqref{eq:Ndb} implies that
\[
\sum_{\substack{d\mid\lambda(b)\\p_b\mid d}}\frac{N(d,b)}{d}\le
\frac1{p_b}\sum_{d\mid\lambda(b)}N(d,b)\le\frac{b}{B(b)}.
\]
Suppose now $p_b\nmid d$.
Since $p_b\mid\lambda(b)\mid\varphi(b)$,
we have either $p_b^2\mid b$ or one or more primes $q\equiv1\pmod{p_b}$
divide $b$.  In either case the number of residues mod $b$ with order
not divisible by $p_b$ is at most $\varphi(b)/p_b$.  (Actually, since $\gcd(d,b)=1$,
the case $p_b^2\mid b$ does not occur.)
Thus, 
\[
\sum_{\substack{d\mid\lambda(d)\\p_b\nmid d}}N(d,b)\le\frac{\varphi(b)}{p_b}
\le\frac{b}{B(b)}.
\]
With the above two displays and
\eqref{eq:sbdens}, $\delta(\Sc_b)\le2/(bB(b))$.
Since the sum of $2/(bB(b))$ converges, the
proof is complete.
\end{proof}
\medskip

\noindent{\bf Remark}.
An immediate corollary of Proposition \ref{prop:ssum} is that
\[
\sum_{n\le x}D(n)\sim c_1x,\quad x\to\infty.
\]
\medskip

\begin{theorem}
\label{theorem:density}
Let
\[
c_0=\lim_{k\to\infty}\delta\Big(\bigcup_{2\le b\le k}\Sc_b\Big).
\]
We have $\delta(\Sc)=c_0$.
\end{theorem}
\begin{proof}
First note that Proposition \ref{prop:ssum} implies that $\bigcup_{2\le b\le k}\Sc_b$
has an asymptotic density, so that $c_0$ exists and $c_0\le 1$.  
For a given integer $b\ge2$, we have seen in the proof of Proposition \ref{prop:ssum}
that $\Sc_b$ is the union of $N(d,b)$
residue classes mod $b^2d$, where $d$ runs over the divisors of $\lambda(b)$
that are coprime to $b$ and $N(d,b)$ is the number of residues mod $b$ of multiplicative
order $d$.   Note that $b^2d<b^3$.  It follows from a complete
inclusion-exclusion argument that the number of $n\le x$ in $\bigcup_{2\le b\le(\log x)^{1/3}}\Sc_b$
is $(c_0+o(1))x$ as $x\to\infty$.  It thus suffices to prove that the number of
$n\le x$ with $n\in \Sc_b$ for some $b>(\log x)^{1/3}$ is $o(x)$ as $x\to\infty$.

  Let $\epsilon(x)\downarrow0$ arbitrarily slowly.
It follows from Erd\H os \cite{E36} that but for $o(x)$ integers $n\le x$,
$n$ has no divisors in the interval $(x^{1/2-\epsilon(x)},x^{1/2+\epsilon(x)})$.
In particular, but for $o(x)$ integers $n\le x$, if $n=ab$ we may
assume that either $a\le x^{1/2}/B(x)$ or $b\le x^{1/2}/B(x)$, where
as before, $B(x)=\exp(\sqrt{\log x})$.

We first consider numbers $n\le x$ with $n\in\Sc_b$ and $(\log x)^{1/3}<b\le x^{1/2}/B(x)$; the argument here is mostly in parallel with the proof of Proposition \ref{prop:ssum}.

Using \cite{BFPS}, the number of integers $b\in(e^j,e^{j+1}]$ with
$P(\lambda(b))\le e^{\sqrt{j+1}}$ is $\ll e^{j-\sqrt{j}}$, so the number of
integers $n\le x$ divisible by one of these $b$'s is $\ll x/e^{\sqrt{j}}$.
Since the sum of $1/e^{\sqrt{j}}$ for $e^{j+1}>(\log x)^{1/3}$ is $o(1)$ as $x\to\infty$,
there are at most $o(x)$ integers $n\le x$ divisible by some
$b\in((\log x)^{1/3},x^{1/2}/B(x)]$ with $P(\lambda(b))\le B(b)$.


Let $p_b=P(\lambda(b))$ and assume that $p_b>B(b)$.  
Let $d\mid\lambda(b)$ with $\gcd(d,b)=1$ and let $r$ be one of the $N(d,b)$
residue classes mod $bd$ where $l_b(r)=d$ and $br\equiv1\pmod d$.
The number of integers $n=ab\le x$ where $a\equiv r\pmod{bd}$ is
at most $1+x/(b^2d)$, so the number of integers $n=ab\le x$ with
$l_b(a)=d$ and $n\in\Sc_b$ is at most $N(d,b)+xN(d,b)/(b^2d)$.  Using
\eqref{eq:Ndb}, we have
\begin{equation}
\label{eq:Sbcount}
\sum_{\substack{n\le x\\n\in\Sc_b}}1\le b+x\sum_{d\mid\lambda(b)}\frac{N(d,b)}{b^2d}.
\end{equation} 
Since the sum of $b$ for $b\le x^{1/2}/B(x)=o(x)$, we wish to show that
\begin{equation}
\label{eq:bgoal}
\sum_{(\log x)^{1/3}<b\le x^{1/2}/B(x)}\sum_{d\mid\lambda(b)}\frac{N(d,b)}{b^2d}=o(1),\quad x\to\infty.
\end{equation}

By \eqref{eq:Ndb} the  inner sum in \eqref{eq:bgoal} when
$p_b\mid d$ is $\le 1/(bp_b)\le 1/(bB(b))$.  Summing this for $b>(\log x)^{1/3}$
is $o(1)$ as $x\to\infty$.  

Now consider the case $p_b\,\nmid\, d$.  As we have seen in the proof of
Proposition \ref{prop:ssum}, we have 
\[
\sum_{\substack{d\mid\lambda(b)\\p_b\,\nmid\,d}}N(d,b)\le\frac{\varphi(b)}{p_b}.
\]
Thus, the inner sum in \eqref{eq:bgoal} is $\le 1/(bp_b)\le1/(bB(p))$.
Summing on $b>(\log x)^{1/3}$  this is $o(1)$ as $x\to\infty$.

We have just shown that the number of integers $n\le x$ of the form $ab$
where $n\in\Sc_b$ and $(\log x)^{1/3}<b\le x^{1/2}/B(x)$ is $o(x)$ as $x\to\infty$.
It remains to consider the case $a\le x^{1/2}/B(x)$.

The number of integers $n\le x$ of the form $ab$ with $a\le x^{1/2}/B(x)$
and $P(b)\le B(x)$ is
\[
\ll\sum_{a\le x^{1/2}/B(x)}\frac{x}{aB(x)}=o(x), \quad x\to\infty,
\]
using standard estimates on the distribution of smooth numbers
(or even using \cite{BFPS}).
Now say $n\le x$ is of the form $ab$ with $1<a\le x^{1/2}/B(x)$
and $n\in\Sc_b$.  This implies that $a^{ab-1}\equiv1\pmod b$.  Let $q=P(b)$,
which we may assume is $>B(x)$ and note that $l_a(q)\mid ab-1$.
Write $b=qm$ and since $b\equiv m\pmod{q-1}$, we have $l_a(q)\mid am-1$.
We distinguish two cases: $m\le B(x)^{1/2}$, $m>B(x)^{1/2}$.

Suppose that $m\le B(x)^{1/2}$.  
Since  $l_a(q)\mid am-1$, we have
$q\mid a^{am-1}-1$.  For a given choice of $a,m$,
the number of primes $q$ with this property is $\ll am\log a$.
Summing this expression over $a,m$ we get $\ll (x\log x)/B(x)$,
and so the number of integers $ab$ is $o(x)$.

Next suppose that $m>B(x)^{1/2}$, so that $q<x/(aB(x)^{1/2})$.  For $a,q$ given,
the number of $m$ is at most 
$1+x/(aql_a(q))$.  The sum of ``1" over $q$ is no
problem, it is at most $\pi(x/(aB(x)^{1/2})$, and so summing on $a$,
we get $\ll x/B(x)^{1/2}=o(x)$.  If $l_a(q)>B(x)^{1/3}$,
then summing $x/(aql_a(q))<x/(aqB(x)^{1/3})$ is also no problem.
So, suppose that $l_a(q)\le B(x)^{1/3}$.  Since there are at most $k\log a$
primes dividing $a^k-1$, by summing on $k\le B(x)^{1/3}$ we see
that the number of choices for $q$ is at most
$B(x)^{2/3}\log x$.  Since $q>B(x)$, we have the sum of $x/(aq)$
over these $q$'s at most $(x\log x)/(aB(x)^{1/3})$, which is negligible
when summed over $a$.
\end{proof}

\section{Computation of $c_0$ and $c_1$}
\label{seccomp}
 
 We immediately have $0<c_0\le c_1$.  Indeed, the second inequality is
 obvious from the definitions, and the first inequality follows
since $\Sc$ contains all numbers $n>2$ with $n\equiv2\pmod 4$.  Further, it is not hard to get larger lower bounds for $c_0$ via an inclusion-exclusion to
find the density of $\bigcup_{2\le b\le k}\Sc_b$ for small values of $k$.  Doing
this with $k=10$ gives $880651/1260^2\approx0.554706$.

It is somewhat easier to get lower bounds for $c_1$.  We have computed the
sum of $\delta(\Sc_b)$ for $2\le b\le 10^4$, getting $\approx0.934328$.

However, getting numerical upper bounds for $c_0,c_1$ is a challenge.
Below is a table of counts of $\Sc$ up to various powers of 10, the
counts to $10^8$ confirm those of Eldar.  In addition, we report on
the sum of $D(n)$ to various powers of 10.

\begin{footnotesize}
\begin{table}[ht]
\caption{Count of members of $\Sc$ below various bounds and partial sums of $D(n)$.}
\label{Ta:Scount}
\begin{tabular}{|l|l|l|} \hline
Bound & Count & Sum \\ \hline
 10 &  2 & 2\\
 $10^2$ &   52  & 61\\
  $10^3$ & 591 & 822 \\
 $10^4$ & 6169 & 8962 \\
 $10^5$ & 62389 & 92383 \\
 $10^6$ & 625941 & 932490 \\
 $10^7$ & 6265910 & 9352861 \\
 $10^8$ & 62677099 & 93613688\\
 $10^9$ & 626836390 & 936403866$~$\\
 $10^{10}$ & 6268593131 & \\ \hline
   \end{tabular}
\end{table}
\end{footnotesize}

Thus, it may be that $c_0<0.627$ and $c_1<0.937$.  We can at least
rigorously prove that $c_0<1$.  Further numerical evidence is given
at the end of this section.

For a finite abelian group $G$ consider the function $N(G)$ 
defined as follows:
\[
N(G)=\sum_{d\mid\#G}\frac{N(d,G)}{d},\quad\hbox{where }
N(d,G)=\#\{g\in G: g \hbox{ has order }d\}.
\]
Writing $G=G_{p_1}\times\cdots\times G_{p_k}$, where $G_p$ is a
$p$-group and $p_1,\dots,p_k$ are the distinct primes dividing $\#G$, we have
\[
N(G)=\prod_{p\mid\# G}N(G_p).
\]
So to get a formula or inequality for $N(G)$ it suffices to do so in the special
case of a finite abelian $p$-group.  The literature has papers on 
counting cyclic subgroups, which is essentially the same problem.  For example,
see T\'oth \cite{T}.  However, it is not hard to directly prove (see below) the inequality
\begin{equation}
\label{eq:ng}
N(G)\le \frac{\tau(\lambda(G))\# G}{\lambda(G)},
\end{equation}
where $\tau(n)$ is the number of divisors of $n$ and $\lambda(G)$ is the universal exponent for $G$.  In the case of interest
for Ordowski's conjecture, this assertion is
\begin{equation}
\label{eq:ndbsum}
\sum_{d\mid\lambda(b)}\frac{N(d,b)}d\le\frac{\tau(\lambda(b))\varphi(b)}{\lambda(b)}.
\end{equation}
This supplies an alternate approach to proving Proposition \ref{prop:ssum}.

Indeed, we know from \cite{EPS} that there is a positive constant $c$ such
that for all large $n$, $\lambda(n)>(\log n)^{c\log\log\log n}$.  Since 
$\tau(k)<k^{1/2}$ for all large $k$, \eqref{eq:ndbsum} implies that
\[
\sum_{b\ge b_0}\sum_{d\mid\lambda(b)}\frac{N(d,b)}{b^2d}
\le \sum_{b\ge b_0}\frac{\tau(\lambda(b))\varphi(b)}{b^2\lambda(b)}
\le\sum_{b\ge b_0}\frac1{b(\log b)^{\frac c2\log\log\log b}}
\]
for $b_0$ sufficiently large.  This implies the sum in Proposition \ref{prop:ssum}
converges.

Now we show that $c_0<1$.  Indeed, it follows
from the above paragraph that for $b$ sufficiently large, we have
\begin{equation}
\label{eq:c0proof}
\frac{\tau(\lambda(b))}{\lambda(b)}<\frac1{(\log b)^3}.
\end{equation}
Note that for any $k$,
\[
\mathcal{A}_k:=\mathbb{N}\setminus\Big(\bigcup_{2\le b\le k}\Sc_b\Big)
\]
contains all $n$ with least prime factor exceeding $k$, so that
$\delta(\mathcal{A}_k)>1/(2\log k)$ for $k$ sufficiently large.  Thus,
for $k$ sufficiently large,
\begin{align*}
\delta(\mathbb{N}\setminus\Sc)&\ge\delta(\mathcal{A}_k)-\delta(\bigcup_{b>k}\Sc_b)\\
&>\frac1{2\log k}-\sum_{b>k}\delta(\Sc_b)
>\frac1{2\log k}-\sum_{b>k}\frac1{b(\log b)^3},
\end{align*}
using \eqref{eq:ndbsum} and \eqref{eq:c0proof}.  Now
\[
\sum_{b>k}\frac1{b(\log b)^3}<\int_k^\infty\frac1{t(\log t)^3}\,dt=\frac1{2(\log k)^2},
\]
so that for $n$ large, $\delta(\mathbb{N}\setminus\Sc)>1/(3\log k)>0$.
This shows that $c_0<1$ as claimed.

This argument could be used in principle to get a numerical upper bound for $c_0$
that is $<1$, but it likely would not be a very good bound.  

Here is a proof of  \eqref{eq:ng}.
\begin{lemma}
\label{lem}
Let $G$ be a finite abelian $p$ group of order $p^n$ and with exponent $p^\lambda$.
Then for $0\le j\le\lambda$, $N(p^j,G)/p^j\le p^{n-\lambda}$, and
$N(G)\le \tau(p^\lambda)p^{n-\lambda}$.
\end{lemma}
\begin{proof}
The second assertion clearly follows from the first one, since 
$\tau(p^\lambda)=\lambda+1$.
So, we concentrate on the first assertion, which we prove by induction.
Write $G$ as $C_{p^{\lambda_1}}\times\cdots\times C_{p^{\lambda_k}}$,
where $1\le\lambda_1\le\dots\le\lambda_k$, $n=\lambda_1+\cdots+\lambda_k$,
and $\lambda=\lambda_k$.   For our base cases we have $j=0$ or $k=1$,
the lemma being clear in either case.   
Now assume
that $j\le\lambda_1$.  Then $N(p^j,G)=p^{jk}-p^{(j-1)k}<p^{jk}$,
so that $N(p^j,G)/p^j<p^{j(k-1)}$.  Now
\[
j(k-1)\le\lambda_1(k-1)\le \lambda_1+\cdots+\lambda_{k-1}=n-\lambda_k,
\]
so the lemma holds in this case.  

We assume the lemma holds for $p$-groups of order smaller than $p^n$.
Suppose $G$ has order $p^n$, exponent $p^\lambda$, rank $k\ge2$, and
 assume that  $\lambda\ge j>\lambda_1$.  
Let $G'$ be the same as $G$ except that $C_{p^{\lambda_1}}$ is replaced
with $C_{p^{\lambda_1-1}}$ and let $G''=C_{p^{\lambda_2}}\times\cdots\times C_{p^{\lambda_k}}$.  
An element of order $p^j$ in $G$ is uniquely expressible as $(u,v)$ where
$u$ is an arbitrary element of $C_{p^{\lambda_1}}$ and $v$ is an element
in $G''$ of order $p^j$.  The same goes for $G'$, except $u$ is only
roaming over $C_{p^{\lambda_1-1}}$ instead of $C_{p^{\lambda_1}}$.  Thus,
we have $N(p^j,G)=pN(p^j,G')$.  By the 
induction hypothesis, we have $N(p^j,G')/p^j\le p^{n-1-\lambda}$.
Multiplying both sides by $p$, we have $N(p^j,G)/p^j\le p^{n-\lambda}$,
which completes the proof.
\end{proof}
 
These thoughts ignore the condition that $\gcd(d,b)=1$, but 
it is not hard to remove the local factors corresponding
to primes dividing $\gcd(\lambda(b),b)$.  In particular if $\varphi_0(b)$
is the largest divisor of $\varphi(b)$ that is coprime to $b$ and
$\lambda_0(b)$ is the largest divisor of $\lambda(b)$ coprime to $b$,
then \eqref{eq:ndbsum} can be improved to
\begin{equation}
\label{eq:reduced}
\delta(\Sc_b)\le \frac{\tau(\lambda_0(b))\varphi_0(b)}{\lambda_0(b)b^2}.
\end{equation}

We have summed this bound for all $b$ with $10^4<b\le 10^6$, getting
$\approx0.00638378$, with $10^4<b\le10^7$,  getting
$\approx0.00673006$, and with $10^4<b\le 2\cdot10^7$, getting
$\approx0.00677103$.  It seems reasonable to assume that the infinite
sum of this bound for all $b>10^4$ is $<0.007$.  Assuming this is so,
our rigorous lower estimate of
$0.934328$ for $c_1$ should be within $0.007$ of the true value,
which is indeed consistent with the evidence afforded by our partial
sums of $D(n)$ in Table \ref{Ta:Scount}.  

One can also try to use these methods to get a numerical estimation for
$c_0$, however, the rigorous estimation from below is difficult.  As
mentioned above, the density of $\bigcup_{2\le b\le 10}\Sc_b$ is about
$0.554706$.  To get a reasonable bound one would want to at least
replace ``10" with ``100" here.  In estimating the tail one can ignore
{\it imprimitive} values of $b$, namely a value of $b$ with $\Sc_b\subset\Sc_{b_0}$
for some $2\le b_0<b$.  For example, if $b=a_0b_0$ where $a_0,b_0\ge2$
and $a_0\equiv1\pmod{b_0}$, then $\Sc_b\subset\Sc_{b_0}$.  In particular,
this holds whenever $b\equiv2\pmod 4$ with $b>2$, or when $b=3m$
with $m>1$ and $m\equiv1\pmod3$.

We have shown that the density of the set $\Sc$ of $n$ with $D(n)\ge1$ exists.
We mention that our results show that the set of numbers $n$ with
$D(n)\ge k$, for any fixed $k$, has a positive asymptotic density.
To see this, note that if $n\equiv p\pmod{p^2}$ for each
of the first $k$ primes (or any set of $k$ primes), then $D(n)\ge k$.
A complicated inclusion-exclusion shows that the density exists.

\section{Pseudoprimes in residue classes}
\label{secres}

We begin with a proof of necessity of the condition from the Introduction
for a residue class to contain a base-$a$ pseudoprime.  Recall the notation
$n_a$ as the largest divisor of $n$ that is coprime to $a$.
\begin{lemma}
\label{lem:resclass}
Suppose $a,r,m$ are integers with $a\ge2$ and $m>0$.  Let
$g=\gcd(r,m)$ and $h=\gcd(l_a(g_a),m)$.  If there is an integer
$n\equiv r\pmod m$ with $n$ a base-$a$ pseudoprime, then $h\mid r-1$,
$g/g_a\mid a$, and in the case that $g$ is even,
there is an integer $k\equiv r\pmod m$ with the
Jacobi symbol $(a/k_{2a})=1$.
\end{lemma}
\begin{proof}
Suppose $n\equiv r\pmod m$ is a base-$a$ pseudoprime.  Then
$a^n\equiv a\pmod n$ and this implies that $a^n\equiv a\pmod{g_a}$.
Since $\gcd(a,g_a)$ is 1, we thus have $a^{n-1}\equiv 1\pmod{g_a}$.
Thus, $l_a(g_a)\mid n-1$, so that $h\mid n-1$.  
We have $n-1\equiv r-1\pmod m$, so
that $n-1\equiv r-1\pmod h$, which implies $r-1\equiv0\pmod h$.  Also,
write $g$ as $ug_a$, so that $u=g/g_a$ is the largest divisor of $a$
all of whose prime factors also divide $a$.  Since $u\mid n$, the congruence $a^n\equiv a\pmod n$ implies that $u\mid a$ (if some prime divides $u$
to a higher exponent than it divides $a$, then $a^n$ and $n$ both have
more factors of this prime than does $a$, a contradiction).
This
proves the first part of the condition.  Now suppose that $g$ is even.
Then $n$ is even, so that $a^n$ is a square mod $n$.  Since $a^n\equiv a\pmod n$,
we have that $a$ too is a square mod $n$.  In particular $a$ is a square
modulo the largest odd divisor of $n$ coprime to $a$, namely $n_{2a}$.
Thus, $(a/n_{2a})=1$.  This completes the proof.
\end{proof}

As mentioned in the Introduction, we conjecture that the conditions of
Lemma \ref{lem:resclass} are not only sufficient for there to be a pseudoprime
base $a$ in the residue class $r\pmod m$, but sufficient for there to be infinitely
many.  We conjecture this based not only on the fact that
it has been proved in many cases, but on the Erd\H os heuristic in
\cite{E}. 

Let us illustrate this heuristic in the case of (base-2) pseudoprimes $n\equiv0\pmod 2$.  We already know that there are infinitely many, but the Erd\H os heuristic
implies the number of them up to $x$ is $>x^{1-\epsilon}$ for any fixed
$\epsilon>0$ and $x$ sufficiently large depending on $\epsilon$.  Consider
primes $p\le y$ with $P(p-1)<y^\epsilon$ and $p\equiv7\pmod8$.
Without the congruence condition it is already conjectured that this
entails a positive proportion of the primes to $y$, just as we know unconditionally
that there is a positive proportion of integers $n\le y$ with $P(n)<y^\epsilon$.
Adding in the congruence condition mod 8 for primes should not matter, and it 
provably doesn't matter when counting integers.  So, assume there are
at least $c_\epsilon\pi(y)$ primes $p\le y$ with $P(p-1)<y^\epsilon$ and
$p\equiv7\pmod 8$, where $c_\epsilon>0$ and $y$ is sufficiently
large depending on the choice of $\epsilon$.  Say the set of primes is
${\mathcal P}_\epsilon(y)$.

Let $x=y^y$ and take subsets of
$\mathcal{P}_\epsilon(y^{1/\epsilon})$ of size $\lfloor\epsilon \log(x/2)/\log y\rfloor$.  
Multiply the primes in each subset, so in this way, each such subset 
corresponds to an integer $n\le x/2$.  Since we are assuming that
$\#\mathcal{P}_\epsilon(y^{1/\epsilon})\ge c_\epsilon\pi(y^{1/\epsilon})$,
the number of subsets formed in this
way is $x^{1-\epsilon+o(1)}$.  Is $2n$ a pseudoprime?  For this to be so
we would need $l(n)\mid 2n-1$, that is, $2n\equiv 1\pmod{l(n)}$.  This
condition forces $l(n)$ to be odd, but at least we already know this
since the primes dividing $n$ are all $\equiv7\pmod 8$, which 
implies that $l(n)\mid\lambda(n)/2$ and that $\lambda(n)\equiv2\pmod4$.  
Let $L$ be
the lcm of all prime powers $p^a\le y^{1/\epsilon}$ with $2<p<y$.
Then $L <(y^{1/\epsilon})^{\pi(y)}=x^{o(1)}$.  The ``probability" that
$2n\equiv1\pmod L$ should be about $1/L$.  Assuming this, the
``expected" number of pseudoprimes constructed this way is at least
$x^{1-\epsilon+o(1)}$.

\medskip

Tables \ref{Ta:psp2inclass8} to \ref{Ta:psp2inclass20}
shows the counts
of pseudoprimes to base 2 for even moduli up to 20.
Compare with the first columns of Table 4 in \cite{PSW}.
The ``Fraction'' column gives the fraction of pseudoprimes
in that class below $10^{16}$.
The symbol ``-'' means that there are no pseudoprimes in that class
due to the conditions in Lemma \ref{lem:resclass}.
The symbol ``na'' in the last column means that that count is not available.

Observe that the odd pseudoprimes far outnumber the even ones
for numbers of the sizes we can compute.  It would be nice to prove
that this continues to hold as one counts to higher levels.

As we mentioned in \cite{PSW}, for most moduli $m$,
the residue class $1\pmod m$ is most popular.
In that work, which gave the counts up to $25\cdot10^9$,
we said that the first exception was $m=37$, which had
more pseudoprimes in class 0 than in class $1\pmod{37}$.
Additional computing reported here finds that $1\pmod{37}$
had more pseudoprimes than $0\pmod{37}$ already at $10^{14}$.

Important table entries for judging the conditions in Lemma \ref{lem:resclass}
are the zero counts.  We list the first few classes with no pseudoprimes
up to $10^{16}$ in Table \ref{Ta:zero1}.

The reader can check using Lemma \ref{lem:resclass} that the residue
classes in Table \ref{Ta:zero1} contain no pseudoprime to base 2.
We searched all moduli $m\le 300$ for empty residue classes up to
$10^{16}$ and found only those predicted by Lemma \ref{lem:resclass}.

Feitsma has computed all odd pseudoprimes to base 2 below $2^{64}$.
They are available at the url \cite{Feitsma}.
We computed the even pseudoprimes to base 2 below $10^{16}$ on two
compute clusters at Purdue University.
The algorithm tested the congruence $2^n\equiv2\pmod n$ for every
$n\equiv2$ or $14\pmod{16}$, except for multiples of 9.
It would have run in about half of the time if we had replaced the condition
on multiples of 9 with the condition that $\gcd(n,2145)=1$.  Also, the 
methods of Feitsma (based on earlier work of Galway) could be applied
to the even pseudoprime count, giving further speed-ups.

\begin{footnotesize}
\begin{table}[ht]
\caption{Number of pseudoprimes to base 2 below various limits in residue classes
mod 2, 4, 6, and 8.}
\label{Ta:psp2inclass8}
\begin{tabular}{|rl|r|r|r|r|r|} \hline
Mod & Class & $\le10^8$ & $\le10^{12}$ & $\le10^{16}$ & Fraction & odd $\le2^{64}$ \\ \hline
 2 &  0 &      7 &     155 &    2045 & 0.000431 & na \\
   &  1 &   2057 &  101629 & 4744920 & 0.999569 & 118968378 \\ \hline
 4 &  0 &      - &       - &       - & - & - \\
   &  1 &   1781 &   90317 & 4215953 & 0.888137 & 104532818 \\
   &  2 &      7 &     155 &    2045 & 0.000431 & na \\
   &  3 &    276 &   11312 &  528967 & 0.111433 & 14435560 \\ \hline
 6 &  0 &      - &       - &       - & - & - \\
   &  1 &   1667 &   86672 & 4074420 & 0.858321 & 101153215 \\
   &  2 &      0 &      12 &      72 & 0.000015 & na \\
   &  3 &    117 &    2251 &   44084 & 0.009287 & 532193 \\
   &  4 &      7 &     143 &    1973 & 0.000416 & na \\
   &  5 &    273 &   12706 &  626416 & 0.131961 & 17282970 \\ \hline
 8 &  0 &      - &       - &       - & - & - \\
   &  1 &   1144 &   60415 & 2869324 & 0.604454 & 70734813 \\
   &  2 &      4 &      84 &    1030 & 0.000217 & na \\
   &  3 &    131 &    5646 &  264955 & 0.055816 & 7220309 \\
   &  4 &      - &       - &       - & - & - \\
   &  5 &    637 &   29902 & 1346629 & 0.283682 & 33798005 \\
   &  6 &      3 &      71 &    1015 & 0.000214 & na \\
   &  7 &    145 &    5666 &  264012 & 0.055617 & 7215251 \\ \hline
\end{tabular}
\end{table}
\end{footnotesize}

\begin{footnotesize}
\begin{table}[ht]
\caption{Number of pseudoprimes to base 2 below various limits in residue classes
mod 10, 12, and 14.}
\label{Ta:psp2inclass10}
\begin{tabular}{|rl|r|r|r|r|r|} \hline
Mod & Class & $\le10^8$ & $\le10^{12}$ & $\le10^{16}$ & Fraction & odd $\le2^{64}$ \\ \hline
10 &  0 &      - &       - &       - & - & - \\
   &  1 &   1082 &   61119 & 2969756 & 0.625612 & 73942273 \\
   &  2 &      0 &      14 &     100 & 0.000021 & na \\
   &  3 &    255 &   12198 &  565493 & 0.119127 & 14942850 \\
   &  4 &      0 &      14 &     112 & 0.000024 & na \\
   &  5 &    203 &    5695 &  160728 & 0.033859 & 2517967 \\
   &  6 &      6 &     116 &    1735 & 0.000365 & na \\
   &  7 &    286 &   12643 &  597165 & 0.125799 & 15879976 \\
   &  8 &      1 &      11 &      98 & 0.000021 & na \\
   &  9 &    231 &    9974 &  451778 & 0.095172 & 11685312 \\ \hline
12 &  0 &      - &       - &       - & - & - \\
   &  1 &   1436 &   77269 & 3641316 & 0.767083 & 89412801 \\
   &  2 &      0 &      12 &      72 & 0.000015 & na \\
   &  3 &      6 &      90 &    1048 & 0.000221 & 7743 \\
   &  4 &      - &       - &       - & - & - \\
   &  5 &    234 &   10887 &  531601 & 0.111988 & 14595567 \\
   &  6 &      - &       - &       - & - & - \\
   &  7 &    231 &    9403 &  433104 & 0.091238 & 11740414 \\
   &  8 &      - &       - &       - & - & - \\
   &  9 &    111 &    2161 &   43036 & 0.009066 & 524450 \\
   & 10 &      7 &     143 &    1973 & 0.000416 & na \\
   & 11 &     39 &    1819 &   94815 & 0.019974 & 2687403 \\ \hline
%
14 &  0 &      1 &      28 &     363 & 0.000076 & na \\
   &  1 &    757 &   42605 & 2155951 & 0.454175 & 54972365 \\
   &  2 &      1 &       8 &     119 & 0.000025 & na \\
   &  3 &    230 &   11111 &  510841 & 0.107614 & 13250508 \\
   &  4 &      0 &      12 &     120 & 0.000025 & na \\
   &  5 &    212 &   10315 &  476087 & 0.100293 & 12230634 \\
   &  6 &      2 &      12 &     117 & 0.000025 & na \\
   &  7 &    228 &    8546 &  288424 & 0.060760 & 5156009 \\
   &  8 &      0 &      65 &    1073 & 0.000226 & na \\
   &  9 &    218 &    9407 &  420766 & 0.088639 & 10637121 \\
   & 10 &      2 &      14 &     124 & 0.000026 & na \\
   & 11 &    184 &    9178 &  409825 & 0.086334 & 10310802 \\
   & 12 &      1 &      16 &     129 & 0.000027 & na \\
   & 13 &    228 &   10467 &  483026 & 0.101755 & 12410939 \\ \hline
\end{tabular}
\end{table}
\end{footnotesize}

\begin{footnotesize}
\begin{table}[ht]
\caption{Number of pseudoprimes to base 2 below various limits in residue classes
mod 16 and 18.}
\label{Ta:psp2inclass16}
\begin{tabular}{|rl|r|r|r|r|r|} \hline
Mod & Class & $\le10^8$ & $\le10^{12}$ & $\le10^{16}$ & Fraction & odd $\le2^{64}$ \\ \hline
16 &  0 &      - &       - &       - & - & - \\
   &  1 &    716 &   39177 & 1896100 & 0.399434 & 47068200 \\
   &  2 &      4 &      84 &    1030 & 0.000217 & na \\
   &  3 &     65 &    2795 &  132181 & 0.027845 & 3609796 \\
   &  4 &      - &       - &       - & - & - \\
   &  5 &    320 &   15334 &  696877 & 0.146805 & 17571790 \\
   &  6 &      - &       - &       - & - & - \\
   &  7 &     76 &    2901 &  132347 & 0.027880 & 3609439 \\
   &  8 &      - &       - &       - & - & - \\
   &  9 &    428 &   21238 &  973224 & 0.205020 & 23666613 \\
   & 10 &      - &       - &       - & - & - \\
   & 11 &     66 &    2851 &  132774 & 0.027970 & 3610513 \\
   & 12 &      - &       - &       - & - & - \\
   & 13 &    317 &   14568 &  649752 & 0.136877 & 16226215 \\
   & 14 &      3 &      71 &    1015 & 0.000214 & na \\
   & 15 &     69 &    2765 &  131665 & 0.027737 & 3605812 \\ \hline
%
18 &  0 &      - &       - &       - & - & - \\
   &  1 &    990 &   54852 & 2654508 & 0.559201 & 65743806 \\
   &  2 &      0 &       5 &      24 & 0.000005 & na \\
   &  3 &     54 &    1117 &   21926 & 0.004619 & 266159 \\
   &  4 &      1 &      20 &     247 & 0.000052 & na \\
   &  5 &    101 &    4197 &  208745 & 0.043974 & 5762593 \\
   &  6 &      - &       - &       - & - & - \\
   &  7 &    341 &   15987 &  709937 & 0.149556 & 17704708 \\
   &  8 &      0 &       6 &      27 & 0.000006 & na \\
   &  9 &      - &       - &       - & - & - \\
   & 10 &      6 &      98 &    1488 & 0.000313 & na \\
   & 11 &     90 &    4287 &  208982 & 0.044024 & 5760564 \\
   & 12 &      - &       - &       - & - & - \\
   & 13 &    336 &   15833 &  709975 & 0.149564 & 17704701 \\
   & 14 &      0 &       1 &      21 & 0.000004 & na \\
   & 15 &     63 &    1134 &   22158 & 0.004668 & 266034 \\
   & 16 &      0 &      25 &     238 & 0.000050 & na \\
   & 17 &     82 &    4222 &  208689 & 0.043963 & 5759813 \\ \hline
\end{tabular}
\end{table}
\end{footnotesize}

\begin{footnotesize}
\begin{table}[ht]
\caption{Number of pseudoprimes to base 2 below various limits in residue classes
mod 20.}
\label{Ta:psp2inclass20}
\begin{tabular}{|rl|r|r|r|r|r|} \hline
Mod & Class & $\le10^8$ & $\le10^{12}$ & $\le10^{16}$ & Fraction & odd $\le2^{64}$ \\ \hline
20 &  0 &      - &       - &       - & - & - \\
   &  1 &    943 &   55255 & 2711430 & 0.571192 & 67162651 \\
   &  2 &      0 &      14 &     100 & 0.000021 & na \\
   &  3 &     33 &    1558 &   76876 & 0.016195 & 2162054 \\
   &  4 &      - &       - &       - & - & - \\
   &  5 &    203 &    5695 &  160728 & 0.033859 & 2517967 \\
   &  6 &      6 &     116 &    1735 & 0.000365 & na \\
   &  7 &     69 &    2505 &  127520 & 0.026863 & 3630971 \\
   &  8 &      - &       - &       - & - & - \\
   &  9 &    196 &    8589 &  385533 & 0.081217 & 9822399 \\
   & 10 &      - &       - &       - & - & - \\
   & 11 &    139 &    5864 &  258326 & 0.054419 & 6779622 \\
   & 12 &      - &       - &       - & - & - \\
   & 13 &    222 &   10640 &  488617 & 0.102933 & 12780796 \\
   & 14 &      0 &      14 &     112 & 0.000024 & na \\
   & 15 &      - &       - &       - & - & - \\
   & 16 &      - &       - &       - & - & - \\
   & 17 &    217 &   10138 &  469645 & 0.098936 & 12249005 \\
   & 18 &      1 &      11 &      98 & 0.000021 & na \\
   & 19 &     35 &    1385 &   66245 & 0.013955 & 1862913 \\ \hline
\end{tabular}

\vspace{1cm}

\caption{List of residue classes for even moduli up to 26 with no pseudoprimes to base 2 up to $10^{16}$.}
\label{Ta:zero1}
\begin{tabular}{|rl|} \hline
Modulus & Empty classes \\ \hline
4 & 0 \\
6 & 0 \\
8 & 0 4 \\
9 & 0 \\
10 & 0 \\
12 & 0 4 6 8 \\
16 & 0 4 6 8 10 12 \\
18 & 0 6 9 12 \\
20 & 0 4 8 10 12 15 16 \\
21 & 0 14 \\
22 & 0 \\
24 & 0 4 6 8 12 16 18 20 \\
25 & 0 \\
26 & 0 \\ \hline
\end{tabular}
\end{table}
\end{footnotesize}

\FloatBarrier
{\noindent{\bf Dedication.}
Our proof of Ordowski's conjecture
 bears some resemblance to a series of papers of Aleksandar Ivi\'c
 \cite[Ch.\ 6]{DI}, \cite{I}, \cite{IP} dealing with tight estimates for the reciprocal sum of the largest
 prime factor of an integer.  We trust he would have enjoyed the connection,
 and we dedicate this paper to his memory.
 
 In addition to Professor Ivi\'c, the year 2020 saw the passing of too many people.
 Among these were John Conway
 and Richard Guy.  They had a deep interest in pseudoprimes, for example,
 Section A12 of \cite{G} is devoted entirely to this subject.  With Schneeberger
 and Sloane, they had a quite remarkable paper \cite{CG} on pseudoprimes.
 For each integer $a\ge0$, let $n_a$ be the least composite number with
 $a^{n_a}\equiv a\pmod{n_a}$.  Then of course $n_a\le 561$, the first
 Carmichael number.  What they showed is the remarkable fact that
 the sequence $n_0,n_1,\dots$ is periodic with period $23\#277\#$,
 where $p\#$ is the product of the primes up to $p$.  Who knew?}


\end{document}